\documentclass[letterpaper, 10pt, conference]{IEEEconf}
\IEEEoverridecommandlockouts 
\usepackage{amssymb,amsmath,times,graphicx,subfigure}

\newcommand{\norm}[1]{\ensuremath{\left\| #1 \right\|}}

\newcommand{\bracket}[1]{\ensuremath{\left[ #1 \right]}}
\newcommand{\braces}[1]{\ensuremath{\left\{ #1 \right\}}}
\newcommand{\parenth}[1]{\ensuremath{\left( #1 \right)}}

\newcommand{\refeqn}[1]{(\ref{eqn:#1})}
\newcommand{\reffig}[1]{Fig. \ref{fig:#1}}
\newcommand{\tr}[1]{\mbox{tr}\ensuremath{\negthickspace\bracket{#1}}}
\newcommand{\deriv}[2]{\ensuremath{\frac{\partial #1}{\partial #2}}}
\newcommand{\SO}{\ensuremath{\mathrm{SO(3)}}}

\renewcommand{\Re}{\ensuremath{\mathbb{R}}}
\renewcommand{\S}{\ensuremath{\mathbb{S}}}

\title{\LARGE \bf
Global Symplectic Uncertainty Propagation on $\SO$}

\author{Taeyoung Lee\authorrefmark{1}\authorrefmark{2}, Melvin Leok\authorrefmark{1}, and N. Harris McClamroch\authorrefmark{2}%
\thanks{Taeyoung Lee, N. Harris McClamroch, Aerospace Engineering, University of Michigan, Ann Arbor, MI 48109 {\tt \{tylee,nhm\}@umich.edu}}%
\thanks{Melvin Leok, Mathematics, Purdue University, West Lafayette, IN 47907 {\tt mleok@math.purdue.edu}}%
\thanks{\textsuperscript{\footnotesize\ensuremath{*}}This research has been supported in part by NSF under grants DMS-0504747, DMS-0726263, and DMS-0714223.}
\thanks{\textsuperscript{\footnotesize\ensuremath{\dagger}}This research has been supported in part by NSF under grant CMS-0555797.}
}

\begin{document}
\maketitle \thispagestyle{empty} \pagestyle{empty}

\begin{abstract}
This paper introduces a global uncertainty propagation scheme for rigid body dynamics, through a combination of numerical parametric uncertainty techniques, noncommutative harmonic analysis, and geometric numerical integration. This method is distinguished from prior approaches, as it allows one to consider probability densities that are global, and are not supported on only a single coordinate chart on the manifold.  The use of Lie group variational integrators, that are symplectic and stay on the Lie group, as the underlying numerical propagator ensures that the advected probability densities respect the geometric properties of uncertainty propagation in Hamiltonian systems, which arise as consequence of the Gromov nonsqueezing theorem from symplectic geometry. We also describe how the global uncertainty propagation scheme can be applied to the problem of global attitude estimation.
\end{abstract}

\section{Introduction}

A nonlinear uncertainty propagation scheme is developed for the dynamics of a rigid body, viewed as evolving on the special orthogonal group $\SO$. Rigid body dynamics is a Hamiltonian flow on a Lie group, and most current attitude uncertainty propagation schemes~\cite{CraMarChe.JGCD07} do not properly take these characteristics into account. Typically, the attitude is represented by unit quaternions, which is problematic for global uncertainty propagation due to the ambiguity introduced by the double cover of $\SO$ by the three-sphere $\S^3$ of unit quaternions. Furthermore, the dynamics are often simplified to kinematic equations, thereby ignoring uncertainties in the angular velocity. As such, most existing techniques are only valid over time periods when the uncertainties are small.

The method introduced in this paper is focused on developing a global uncertainty propagation method for a rigid body by explicitly considering the characteristics of the Hamiltonian flow and the configuration manifold, without implicitly assuming that the uncertainty is localized, nor that the uncertainty distribution is fully supported in a single coordinate chart on the manifold.

Although it is not widely known in the engineering community, the theories of probability and stochastic processes on manifolds have been studied by theoretical statisticians~\cite{Dia.BK88,Elw.BK82,Eme.BK89}. Earlier works on attitude estimation on $\SO$ include \cite{LoEsh.SIAP79}, where a probability density function is expressed using noncommutative harmonic analysis. This idea of using fourier analysis on manifolds has been applied in~\cite{Hen.AS90,Kim.AS98,ParKimZho.IEEECRA05,WanZhoMas.IEEETC06}, and recently, results in~\cite{LoEsh.SIAP79} are extended to include the effects of process noise and sensor parameters in~\cite{Mar.AAS05}.

The Liouville equation describes the evolution of a probability density function in the absence of external diffusion, and can be viewed as the deterministic analogue of the Fokker-Planck equation. When the flow that is advecting the probability density is Hamiltonian, the Liouville equation reduces to an ordinary differential equation~\cite{SchHsiPar.AAS05}. Thus, a probability density function can be propagated using the flow map of the Hamiltonian system. In~\cite{ValPal.AIAA06}, an attitude estimation scheme is developed by linearizing the attitude dynamics along the mean obtained from the probability density function propagated by this property.  However, nonlinearities of the flow imply that numerical methods for propagating uncertainties using linearization rapidly degrade in performance, unless frequent physical measurements are available \cite{CDC07.est}. It is therefore desirable to construct efficient numerical methods for solving the Liouville equation that describes the evolution of a probability density advected by a prescribed Hamiltonian flow.

In this paper, a global attitude uncertainty propagation scheme is developed in the absence of process noise. This is achieved through a synthesis of numerical parametric uncertainty analysis techniques \cite{Xiu2007, XiuHes2005}, noncommutative harmonic analysis \cite{ChKy2001}, and geometric numerical integration~\cite{HaLuWa2006}. We backpropagate sample points along the Hamiltonian flow of the rigid body dynamics, and we use the advected probabilities to reconstruct the probability density using noncommutative harmonic analysis on $\SO$. This is in contrast to Monte Carlo methods, where the sample trajectories are used to compute statistical properties of the advected density.

The Gromov nonsqueezing theorem \cite{Gr1985} from symplectic geometry places fundamental limits on how the uncertainty of a Hamiltonian system evolves \cite{Scetal2007}. It is therefore essential that symplectic methods be used to propagate individual trajectories when analyzing uncertainty propagation in Hamiltonian systems. We use a geometric numerical integrator, referred to as a Lie group variational integrator, that preserves the symplectic property of the Hamiltonian dynamics and the group structure of the configuration space~\cite{CMA07,CMDA07}. The purpose of this paper is to develop a computational method to propagate uncertainties under a Hamiltonian flow on a Lie group.   Our development is for a specific Hamiltonian system, namely the 3D pendulum, that evolves on the Lie group $\SO$.   We also comment on how this development is applicable to attitude estimation. The methodology proposed in this paper should be distinguished from the numerical methods that compute a specific realization of a stochastic Hamiltonian system~\cite{lazarocami-2007, bourabee-2007,Wan.Phd07}.

\section{Attitude Dynamics of a Rigid Body}

\subsection{3D Pendulum}

The 3D pendulum is  nontrivial example of a Hamiltonian system that evolves on a Lie group; this example is used in the subsequent development.    A rigid 3D pendulum is a rigid body supported by a fixed, frictionless pivot, acted on by uniform gravitational force~\cite{SheSanCha.CDC04}. The supporting pivot allows the pendulum three rotational degrees of freedom.

Two reference frames are introduced. An inertial reference frame has its origin at the pivot; the first two axes lie in the horizontal plane and the third axis is vertical in the direction of gravity. A reference frame fixed to the pendulum body is also introduced.  The origin of this body-fixed frame is also located at the pivot.  The configuration manifold is the special orthogonal group $\SO$,
\begin{align*}
    \SO = \braces{ R\in\Re^{3\times 3}\,|\, R^T R=I_{3\times 3},\, \det{R}=1},
\end{align*}
where the rotation matrix $R\in\SO$ represents the linear transformation from the body-fixed frame to the inertial frame.

The dynamics of the 3D pendulum are given by the Euler rigid body equation that includes the moment due to gravity:
\begin{equation}
J \dot\Omega =J\Omega \times \Omega + mg \rho \times R^T e_3, \label{eq:Jw_dot}
\end{equation}
where the angular velocity in the body-fixed frame is denoted by $\Omega\in\Re^3$, the inertia matrix is denoted by $J\in\Re^{3\times 3}$, and the vector $\rho\in\Re^3$ represents the location of the center of mass in the body-fixed frame. The constants $m$ and $g$ denote the mass of the pendulum and the gravitational acceleration, respectively. The kinematic equation is
\begin{equation}
\dot{R} = R S(\Omega).\label{eq:R_dot}
\end{equation}
For a given vector $x\in\Re^3$, the $3\times 3$ skew-symmetric matrix $S(x)$ is defined so that $S(x)y=x\times y$ for all $y\in\Re^3$.

There are two disjoint equilibria when the direction of
gravity in the body fixed frame is collinear with the vector $\rho$; the hanging equilibrium when $R^Te_3=\rho/\norm{\rho}$, and the inverted equilibrium when $R^Te_3=-\rho/\norm{\rho}$. The 3D pendulum exhibits surprisingly rich and complicated attitude dynamics~\cite{CDC07.est}, and is therefore particularly appropriate for demonstrating the properties of our global attitude uncertainty propagation scheme.

\subsection{Lie Group Variational Integrator}
Lie group variational integrators preserve the group structure without the use of local charts, reprojection, or constraints, they are symplectic and momentum preserving, and they exhibit good energy behavior for an exponentially long time period. The following Lie group variational integrator for the attitude dynamics of a rigid body was presented in~\cite{CMA07,CMDA07}:
\begin{gather}
h S(J\Omega_k+\frac{h}{2} M_k) = F_k J_d - J_dF_k^T,\label{eqn:findf0}\\
R_{k+1} = R_k F_k,\label{eqn:updateR0}\\
J\Omega_{k+1} = F_k^T J\Omega_k +\frac{h}{2} F_k^T M_k
+\frac{h}{2}M_{k+1},\label{eqn:updatew0}
\end{gather}
where $J_d\in\Re^{3\times 3}$ is a nonstandard moment of inertia matrix given by $J_d=\frac{1}{2}\mathrm{tr}\!\bracket{J}I_{3\times 3}-J$, and $F_k\in\SO$ is the relative attitude between integration steps. The moment due to the potential is denoted by $M_k=mg\rho\times R_k^T e_3$. The constant $h\in\Re$ is the integration step size, and the subscript $k$ denotes the $k$-th integration step. This integrator yields a map $(R_k,\Omega_k)\mapsto(R_{k+1},\Omega_{k+1})$ by solving \refeqn{findf0} to obtain $F_k\in\SO$ and substituting it into \refeqn{updateR0} and \refeqn{updatew0} to obtain $R_{k+1}$ and $\Omega_{k+1}$. We use these discrete-time equations of motion to propagate the attitude dynamics.

\section{Global Symplectic Uncertainty Propagation on $\SO$}

Suppose that a probability density function for the attitude and the angular velocity of a rigid body at a time $t_k$ is given as $p_k(R,\Omega): \SO\times\Re^3\rightarrow \Re$. In this section, we develop a computational method to propagate this density along the Hamiltonian attitude dynamics assuming that there is no process noise. We represent the propagated probability density function using noncommutative harmonic analysis. A method for visualizing the attitude uncertainty is also discussed.

\subsection{Symplectic Uncertainty Propagation for a Hamiltonian System}

In~\cite{SchHsiPar.AAS05}, it is shown that the probability density function is preserved along a Hamiltonian flow on Euclidean space. In this subsection, we generalize this result to a general Hamiltonian system evolving on a manifold.

Consider a Hamiltonian system on a $2n$-dimensional symplectic manifold $(Q,\omega)$, where $Q$ is a $2n$-dimensional manifold and $\omega:TQ\times TQ\rightarrow\Re$ is a nondegenerate symplectic two-form on $Q$~\cite{MarRat.BK99}. The Liouville volume form $\mu:(TQ)^{2n}\rightarrow\Re$ is defined as the $n$-fold wedge product of the symplectic two-form with itself, $\mu = \frac{(-1)^{n(n-1)/2}}{n!} \omega \wedge \cdots \wedge \omega \quad(n \text{ times})$. In local coordinates, this corresponds to the usual notion of the volume element in  Euclidean spaces. Let $\mathcal{L}_X \mu$ be the Lie derivative of the volume form $\mu$ along a vector field $X:Q\rightarrow TQ$. The divergence of a vector field $X$ on $Q$ is defined as $\mathcal{L}_X \mu = \mathrm{div}_\mu (X)\, \mu$. Thus, the divergence $\mathrm{div}_\mu (X)$ represents the rate of change of a unit volume along the vector field $X$.

In~\cite{ChiKya.BK00}, it is shown that the Fokker-Planck equation for a dynamic system on a manifold in the absence of diffusion terms can be written as,
\begin{align}
    \deriv{p}{t} + \mathrm{div}_\mu \, (pX)=0.
\end{align}
Note that this has the same structure as the Euler equation for the density of compressible fluids. The existence and uniqueness of the solution of the equations of motion provides a property analogous to mass conservation in fluid dynamics. Since $\mathrm{div}_\mu\, (pX) = p \, \mathrm{div}_\mu\,(X)+\mathcal{L}_X p$, the time derivative of the probability density function is given by
\begin{align}
    \frac{d}{dt} p =\deriv{p}{t} + \mathcal{L}_X p = - p\, \mathrm{div}_\mu \, (X).
\end{align}
If the vector field $X$ is the Hamiltonian vector field on $(Q,\omega)$, the divergence vanishes, $\mathrm{div}_\mu \, (X)=0$ according to the Liouville theorem~\cite{MarRat.BK99}. Therefore, the Fokker-Planck equation for a deterministic Hamiltonian system is represented by the ordinary differential equation
\begin{align}
    \frac{d}{dt} p =0.\label{eqn:FPH}
\end{align}

This states that the probability density function is preserved along a Hamiltonian flow without stochastic diffusion effects. More precisely, \refeqn{FPH} implies that the propagated probability density function at $t_{k+1}$ can be explicitly expressed as a composition of the backward flow map and the given probability density function at $t_k$
\begin{align}
    p_{k+1} (R,\Omega) = p_k ( \mathcal{F}^{-1} (R,\Omega)),\label{eqn:pkp}
\end{align}
where $\mathcal{F}^{k}:\SO\times\Re^3\rightarrow \SO\times\Re^3$ represents the $k$ step discrete flow of the attitude dynamics. We can apply this equation recursively to propagate the probability density function over any time interval.

\subsection{Noncommutative Harmonic Analysis on $\SO$}

Equation \refeqn{pkp} provides a method to compute the probability density function at any time based on the probability density function at some prior time. As the flow is nonlinear, it is inefficient to characterize the density using its moments, since the moment expansion may decay slowly.
Here, we propose a computational scheme that represents the probability density function for the attitude dynamics using noncommutative harmonic analysis on $\SO$.

Noncommutative harmonic analysis is a generalization of Fourier analysis on Euclidean spaces to manifolds~\cite{ChiKya.BK00}. This is particularly useful since it provides a mathematical tool to approximate a probability density function based on samples of that function. More precisely, the propagated probability density function for the attitude dynamics of a rigid body can be expressed as
\begin{align}
    p (R,\Omega) = \sum_{l=0}^\infty \frac{2l+1}{(2\pi)^3} \int_{\Re^3} \exp(j\theta\cdot\Omega)\, \tr{P^l(\theta) U^l(R)}\, d\theta,\label{eqn:rec}
\end{align}
where the vector $\theta\in\Re^3$ and non-negative integer $l\in\mathbb{N}\bigcup\{0\}$ are Fourier parameters, and the set of complex matrices $\{P^l(\theta)\in\mathbb{C}^{(2l+1)\times (2l+1)}\}_{l=0}^{\infty}$ is the Fourier spectrum of the density $p(R,\Omega)$. We denote the $l$-th irreducible unitary representation matrix of $\SO$ by $U^l(R)\in\mathbb{C}^{(2l+1)\times (2l+1)}$.
A representation of a group is a homomorphism from the group to the set of invertible matrices, and by the Peter-Weyl theorem~\cite{PetWey.MA27}, the irreducible unitary representations form a complete orthonormal basis for the set of square-integrable functions on the group. The irreducible unitary representations can be expressed in various ways. For example, it can be expressed in terms of the \textit{wigner}-$d$ functions using 3-1-3 Euler angles $\alpha,\beta,\gamma$ as
\begin{align}
    U^l_{m,n}(R(\alpha,\beta,\gamma)) = i^{m-n} e^{-i(m\alpha+n\gamma)} d^l_{mn}(\cos\beta)
\end{align}
for $-l\leq m,n\leq l$~\cite{BieLou.BK81}. A few \textit{wigner}-$d$ functions are given by
\begin{align*}
    d^0(\cos\beta) & = 1,\\
    d^1(\cos\beta) & = \begin{bmatrix}
        \frac{1+\cos\beta}{2} & -\frac{\sin\beta}{\sqrt{2}} & \frac{1-\cos\beta}{2}\\
        \frac{\sin\beta}{\sqrt{2}} & \cos\beta & -\frac{\sin\beta}{\sqrt{2}}\\
        \frac{1-\cos\beta}{2} & \frac{\sin\beta}{\sqrt{2}} & \frac{1+\cos\beta}{2}\end{bmatrix}.
\end{align*}
Higher order \textit{wigner}-$d$ functions can be obtained using a recursive formula~\cite{ChiKya.BK00}.

From the orthnormal property of the irreducible unitary representation, the Fourier spectrum is computed by
\begin{align}
    P^l(\theta) = \int_{\SO} \int_{\Re^3} p(R,\Omega) \exp(-j\theta\cdot\Omega) U^l(R^{-1}) d\Omega\,dR.\label{eqn:fspectrum}
\end{align}
The volume element for the rotation matrix $dR$ represents the Haar measure of $\SO$; it can be written in terms of 3-1-3 Euler angles $\alpha,\beta,\gamma$ as $dR(\alpha,\beta,\gamma) = \frac{1}{8\pi^2} \sin\beta\,d\alpha d\beta d\gamma$.

Substituting \refeqn{pkp}, we can compute the Fourier spectrum of the propagated probability density function. The propagated distribution can be reconstructed by \refeqn{rec}. This is a global particle-based method to construct the propagated probability density function on $\SO\times\Re^3$.

\subsection{Visualization of the Attitude Uncertainty}
Let $p_R:\SO\rightarrow\Re$ be a probability density function on $\SO$. For example, it can be obtained by integrating \refeqn{rec} over $\Re^3$, i.e. $p_R(R)=\int_{\Re^3} p(R,\Omega) d\Omega$. We propose a method for visualizing probability densities on $\SO$ using three copies of two-spheres. The rotation matrix represents a linear transformation from a body fixed frame to an inertial frame. Therefore, the $i$-th column of the rotation matrix $Re_i$ represent the direction of the $i$-th body fixed axis in the inertial frame. Since the vector $Re_i$ lies on the two-sphere $\S^2$, we can visualize uncertainties of $Re_i$ on a sphere either by color shading or by contour lines. Three copies of these spheres, one for each of the body fixed axes, can be used to visualize uncertainties on $\SO$.

We find the marginal probability density of $p_R$ for each column of the rotation matrix. For given $r\in\S^2$, define a subset of $\SO$ as
\begin{align}
    H_i(r) = \braces{ R\in\SO\,\big|\, Re_i = r}.
\end{align}
This is a subgroup of $\SO$ diffeomorphic to $\S^1$, and it represents the set of rotation matrices whose $i$-th column is equal to the given direction $r$. This subgroup can be parameterized by $\theta\in\S^1$. More explicitly, we can find an element $R_i^\circ(r)$ of $H_i(r)$ (for example, we have $R_i^\circ(r)=\exp (\frac{\mathrm{acos} (r\cdot e_i)}{\norm{e_i\times r}} S(e_i\times r))$ if $e_i\times r\neq 0$).  Then, the subgroup $H_i(r)$ is parameterized as
\begin{align}
    H_i (r) = \braces{ R_i^\circ(r) \exp S(\theta e_i)\,\big|\,\theta\in\S}.
\end{align}

The corresponding quotient space is the two-sphere, $\SO/H_i(r)\simeq\S^2$. Using the properties of integration on a quotient space of a Lie group, we have
\begin{align}
    1&=\int_{\SO} p_R(R)\, dR\nonumber\\
     &= \int_{r\in S^2} \parenth{\frac{1}{2\pi}\int_{\theta\in\S^1} p_R(R_i^\circ(r) \exp S(\theta  e_i)\,d\theta}\,dr.
\end{align}
Therefore, the marginal probability density for the $i$-th column of the rotation matrix,  $p_R^i:\S^2\rightarrow\Re$ is given by
\begin{align}
    p_R^i (r) = \frac{1}{2\pi}\int_{\theta\in\S^1} p_R (R_i^\circ(r) \exp (\theta\hat e_i))\,d\theta.\label{eqn:pi}
\end{align}

\begin{figure}
\centerline{
    \subfigure[Density 1]{%
        \includegraphics[width=0.35\columnwidth]{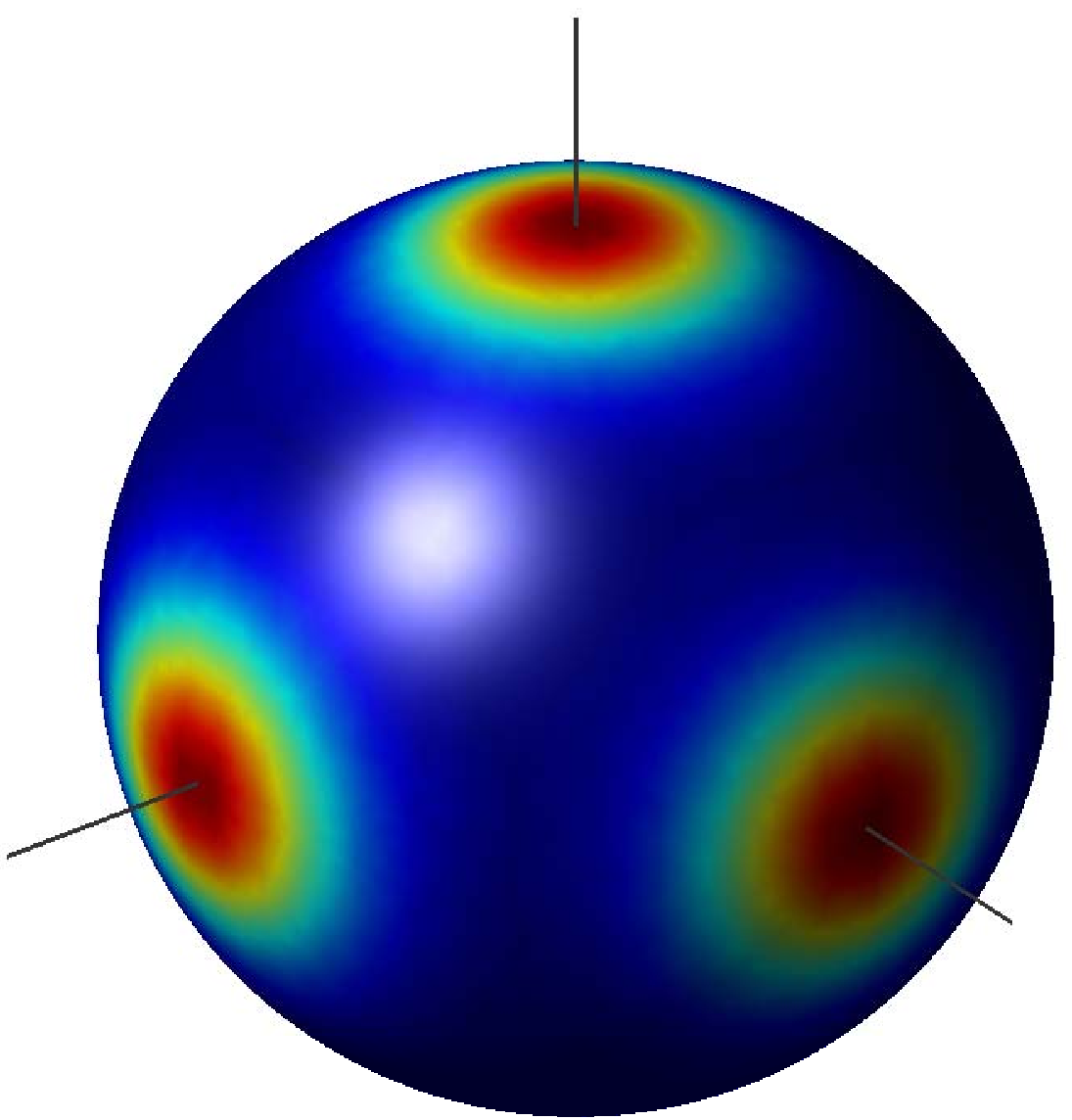}\label{fig:pia}}%
    \hspace*{1.3cm}
    \subfigure[Density 2]{%
        \includegraphics[width=0.35\columnwidth]{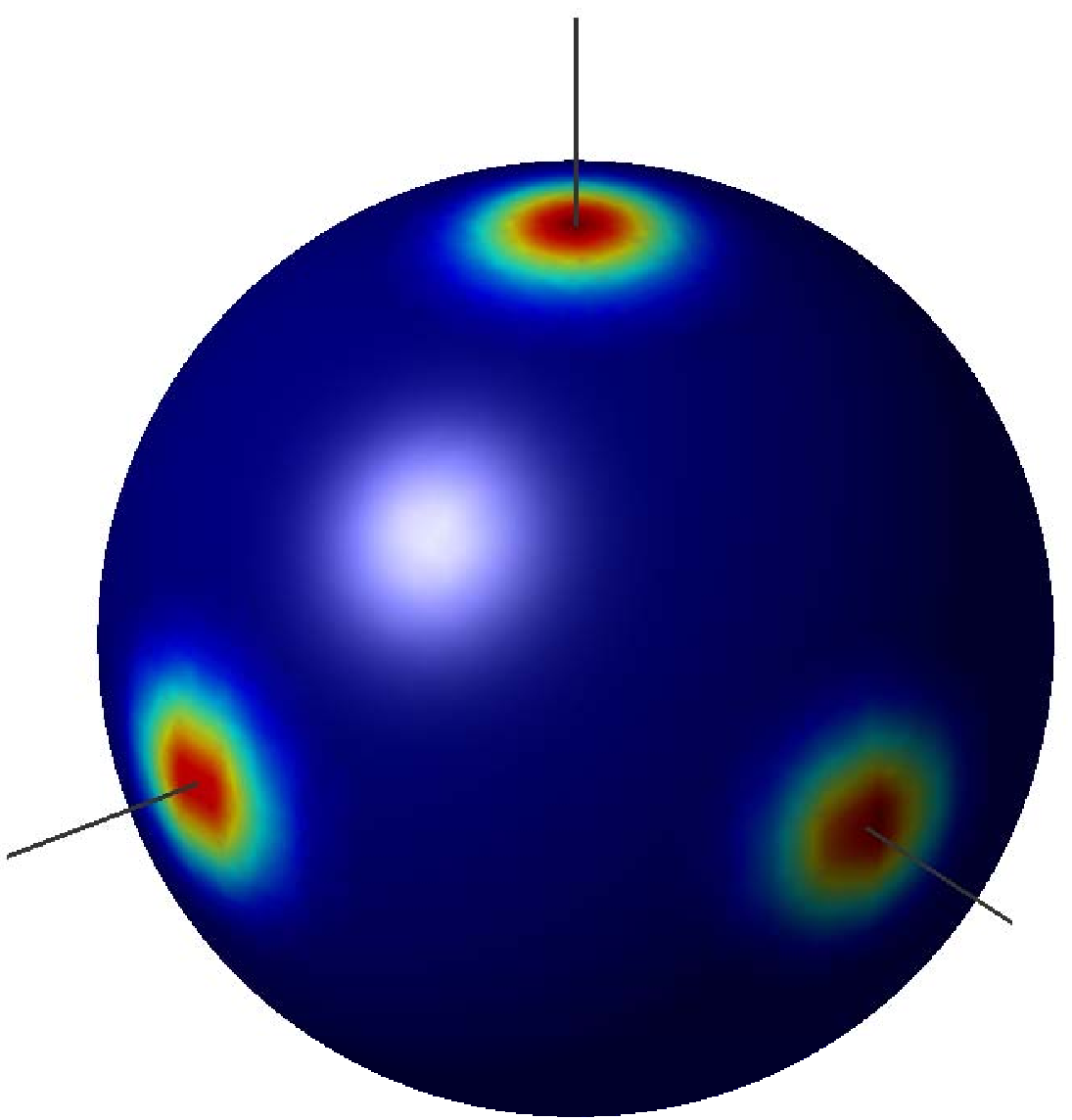}\label{fig:pib}}%
}
\centerline{
    \subfigure[Density 3]{%
        \includegraphics[width=0.35\columnwidth]{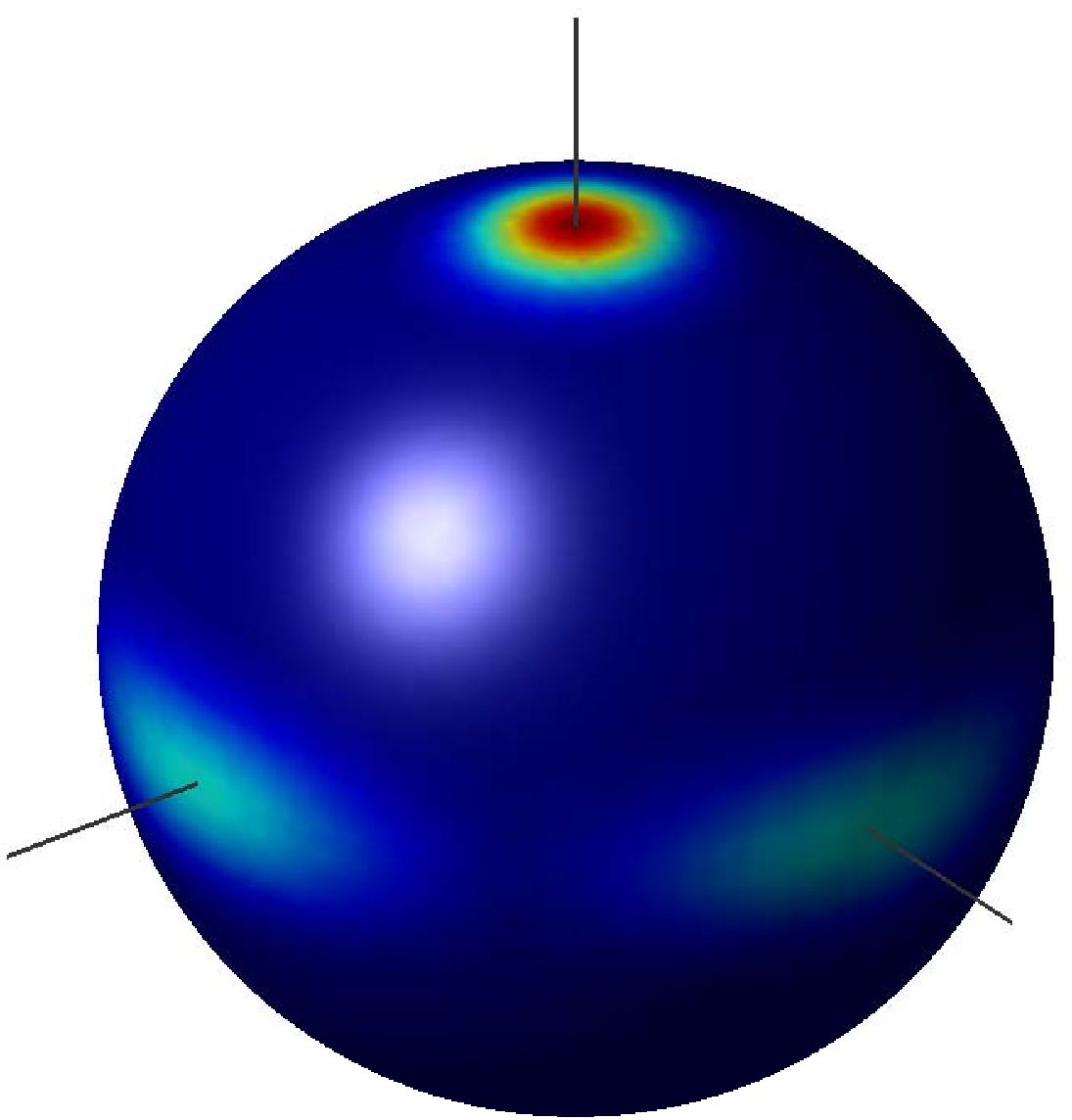}\label{fig:pic}}%
    \hspace*{0.4cm}
    \subfigure[Sample attitudes for the density 3]{%
        \includegraphics[width=0.45\columnwidth]{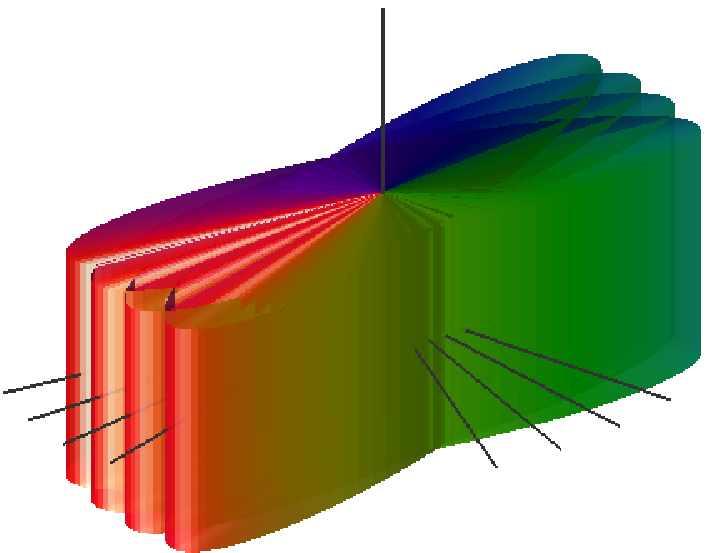}\label{fig:attc}}%
}
\caption{Attitude uncertainty visualization example}\label{fig:vis}
\end{figure}

We plot these marginal probability density functions $p_R^i(r)$, that represent the probability density of the direction of the body fixed axes, on three two-spheres. If the magnitude of uncertainties is small, we can plot uncertainties for each body fixed axis on a single sphere, from which we can intuitively understand the attitude uncertainty of the rigid body. \reffig{vis} shows examples for the attitude probability density visualization. It is easily to see that the second density in \reffig{pib} has smaller variation than the first density in \reffig{pib}. \reffig{attc} shows some sample attitudes that we are likely to obtain from the third density given in \reffig{pic}, which reflect the fact that the direction of the z-axis of the body in the inertial frame is well localized, but there is greater uncertainty in the direction of the x- and y-axes.

\section{Numerical Computations}
In this section, we propagate an initial probability density along the nontrivial dynamics of the 3D pendulum, and we visualize the attitude uncertainty.
The properties of the pendulum are given by
\begin{gather*}
    J=\mathrm{diag}[0.13,0.28,0.17]\,\mathrm{kgm^2},\quad m=1\,\mathrm{kg},\\
    \rho=0.3 e_3\,\mathrm{m},\quad g=9.81\,\mathrm{m/s^2}.
\end{gather*}

The von Mises distribution (also known as the circular normal distribution) is a continuous probability density on the circle, which can be thought of as the circular analogue of the normal density~\cite{Mises.Phy18}.
\begin{align}
    p(\theta) = \frac{1}{2\pi I_0 (\kappa)} \exp \parenth{\kappa \cos(\theta-\overline{\theta})},
\end{align}
where $I_0$ is the zeroth order modified Bessel function of the first kind, given by $I_0(\kappa) = \sum_{i=0}^{\infty} \frac{(1/4\kappa^2)^i}{(\kappa !)^2}$ with parameters $\kappa,\overline\theta$ that determine the shape of the density; as $\kappa$ increases, the density approaches a normal density with mean $\overline\theta$ and variance $\frac{1}{\kappa}$.

For two given rotation matrices $R, \overline R\in\SO$, the quantity $\frac{1}{2}(\tr{\overline R^T R}-1)$ represents the cosine of the rotation angle between the two attitudes. Using this, we define a probability density function on $\SO$ from the von Mises distribution. The probability distribution at the initial time is chosen as
\begin{align}
    p_0(R,\Omega) & = \frac{1}{c} \exp \braces{\frac{1}{2}\kappa\parenth{\tr{\overline R_0^T R}-1}}\nonumber\\
    & \times \exp\braces{-\frac{1}{2} (\Omega-\overline\Omega_0)^T \Sigma^{-1} (\Omega-\overline\Omega_0)},
\end{align}
where $\overline R_0=I_{3\times 3}$, $\overline\Omega_0=[4.14,4.14,4.14]\,\mathrm{rad/s}$, $\Sigma=0.1414^2 I_{3\times 3}$, and $\kappa=8$. The constant $c$ is a scaling factor chosen such that $\int p(R,\Omega)\, dR d\Omega =1$. The corresponding mean $(\overline R_0,\overline\Omega_0)$ yields an irregular, perhaps chaotic, attitude response~\cite{CDC07.est}.

We propagate this initial distribution using \refeqn{pkp} and compute the Fourier spectrum using \refeqn{fspectrum}. The volume integration is approximated by Simpson's rule, and the flow map is computed using the Lie group variational integrator. The application of the Lie group variational integrator is particularly useful since it is symplectic, group structure preserving, and time reversible.

\begin{figure}
\centerline{
    \subfigure[$t=0.0$]{%
        \includegraphics[width=0.4\columnwidth]{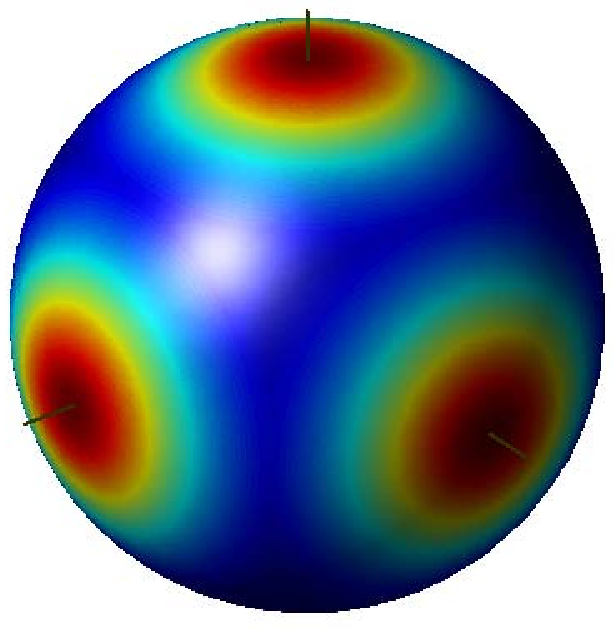}}%
    \hspace*{0.3cm}
    \subfigure[$t=0.1$]{%
        \includegraphics[width=0.4\columnwidth]{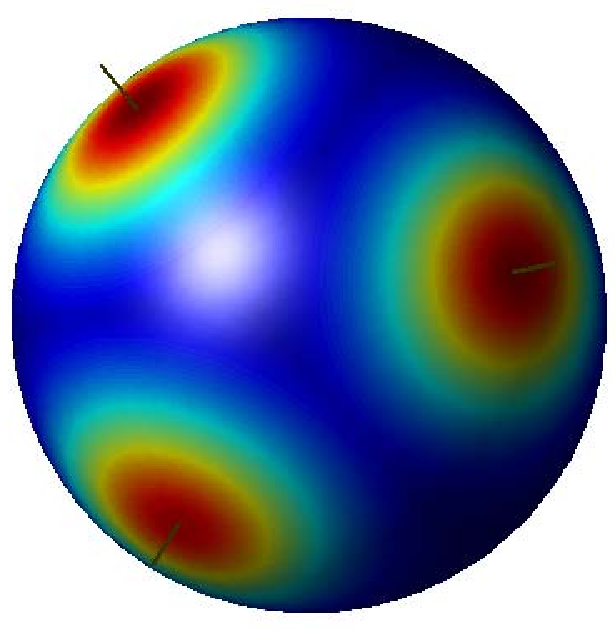}}
}
\centerline{
    \subfigure[$t=0.2$]{%
        \includegraphics[width=0.4\columnwidth]{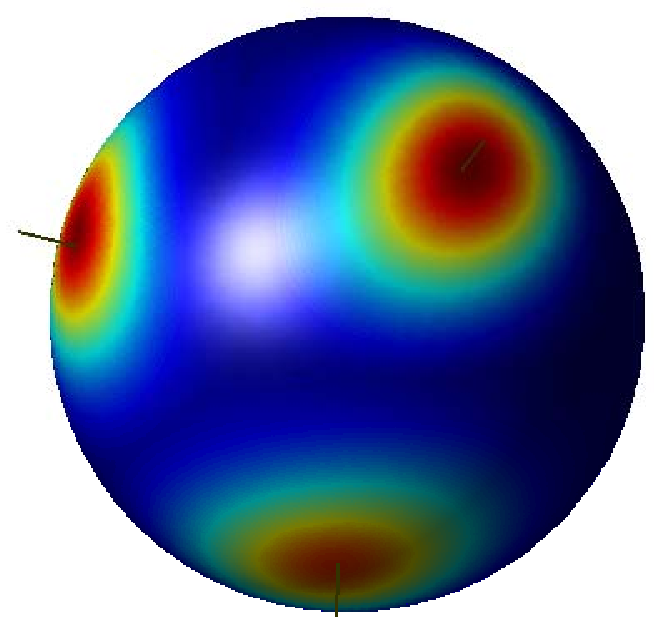}}%
    \hspace*{0.3cm}
    \subfigure[$t=0.4$]{%
        \includegraphics[width=0.4\columnwidth]{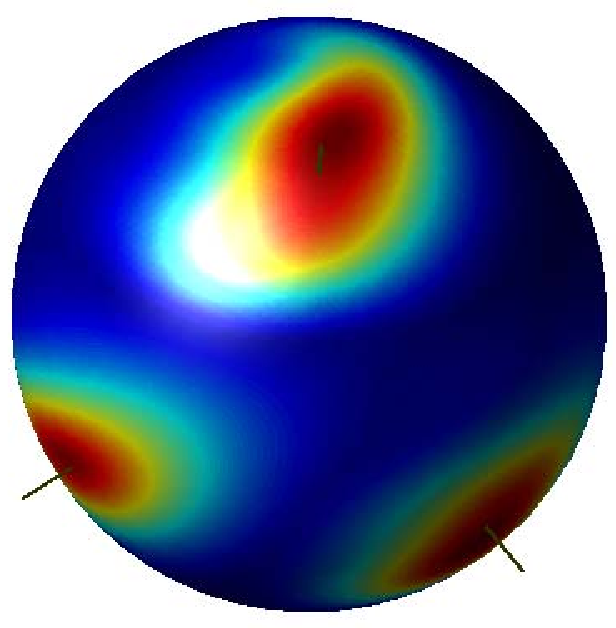}}
}
\centerline{
    \subfigure[$t=1.0$]{%
        \includegraphics[width=0.4\columnwidth]{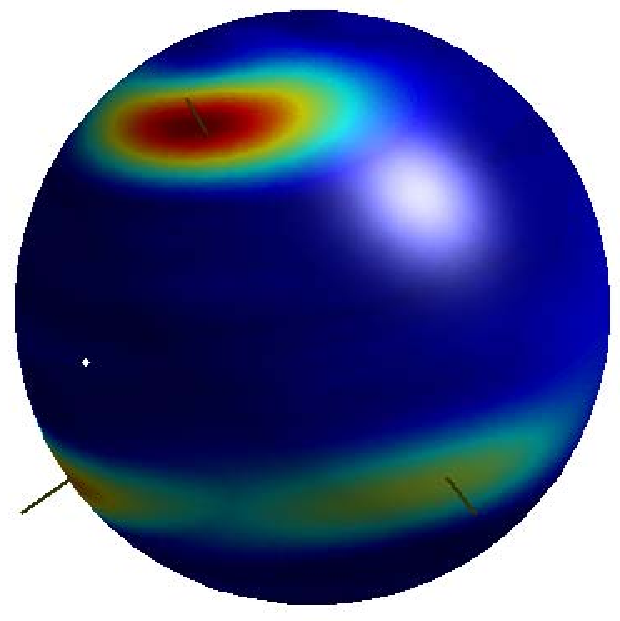}}%
}
\caption{Propagated attitude uncertainty and $\max$ mean attitude}\label{fig:distk}
\end{figure}

We have developed a parallel computing code using the MPI (Message Passing Interface) library, where the domain of integration is divided uniformly and distributed to each processor. This is desirable in terms of minimizing communication between processors and balancing the computational load among the processors. This algorithm has been implemented on 32 AMD Opeteron processors. \reffig{distk} illustrates the propagated attitude uncertainty.

\section{Comments on Global Attitude Estimation}
The presented uncertainty propagation method can be applied to develop an attitude estimation scheme using Bayes rule. We first define a measurement model. We assume that the attitude and the angular velocity are measured by sensors to obtain
\begin{align}
    z_k = H(R_k,\Omega_k)+v_k,\label{eqn:zk}
\end{align}
where $H:\SO\times\Re^3 \rightarrow \Re^m$ is a measurement function, $z_k\in\Re^m$ is the measured value, and $v_k\in\Re^m$ is measurement noise. For example, if we measure a direction to a known object $a\in\S^2$ and the angular velocity, the measurement function can be written as $H(R,\Omega)=[R^T a; \Omega]$. We assume the probability density of the measurement noise $v_k$ is given by $p_{z_k|k}$, and it is an independent process.

The set of all measurements from the initial time to $t_k$ is denoted by the capital letter $Z_k=[z_0,z_1,\ldots,z_k]$. Suppose that we have a probability density function at the $k$-th time conditioned by $Z_k$, i.e. the expression for $p_{k | Z_k}$ is known, and a new measurement is obtained $z_{k+1}$ at $t_{k+1}$. Estimation can be described as finding an expression for $p_{{k+1} | Z_{k+1}}$ given $p_{k | Z_k}$ and $Z_{k+1}=[Z_k, z_{k+1}]$. Using Bayes rule~\cite{Jaz.BK70}, we have
\begin{align}
p_{{k+1}|Z_{k+1}}&(R,\Omega|[Z,z])\nonumber\\
&=\frac{p_{z_{k+1}|{k+1},Z_k}(z|R,\Omega,Z)\,\,p_{{k+1}|Z_k}(R,\Omega|Z)}%
    {p_{z_{k+1}|Z_k}(z|Z)}.\label{eqn:pkpz}
\end{align}
Since the measurement processes are independent, we have $p_{z_{k+1}|{k+1},Z_k}=p_{z_{k+1}|k+1}$, and the propagated density is given by \refeqn{pkp}, i.e. $p_{{k+1}|Z_k}(R,\Omega|Z)=p_{k|Z_k} (\mathcal{F}^{-1} (R,\Omega)|Z)$. The denominator is a normalizing constant that can be computed to satisfy $p_{z_{k+1}|Z_k}(z|Z)=\int_{\SO\times\Re^3} p_{z_{k+1}|k+1}(R,\Omega) p_{k+1|Z_k} (R,\Omega|Z) dRd\Omega$. In summary, the propagated probability density conditioned by the new measurement is given by
\begin{align}
    p_{{k+1}|Z_{k+1}}&(R,\Omega|[Z,z])\nonumber\\
    &=\frac{1}{c}p_{z_{k+1}|{k+1}}(z|R,\Omega)\,\,p_{{k}|Z_k}(\mathcal{F}^{-1}(R,\Omega)|Z)
\end{align}
for a constant $c$. This also can be represented using the harmonic analysis as in \refeqn{rec}.

This expression is the key to solving the attitude estimation problem. This expression for $p_{x_{k+1}|Z_{k+1}}$ contains all the statistical information that can be derived from the attitude dynamics and the available measurements.

\section{Conclusions}

This paper addresses the problem of propagating the attitude uncertainty of a rigid body which evolves on the rotation group. The use of noncommutative harmonic analysis techniques to represent the uncertainty distribution in a global fashion overcomes a fundamental limitation of existing techniques, which implicitly assume that the uncertainty is localized or small. By exploiting the fact that the Liouville equation for a Hamiltonian system reduces to an ordinary differential equation, we are able to adopt a particle based approach for computing the advected probability density, thereby avoiding the computational expense of solving the equation as a numerical partial differential equation. By adopting a Lie group variational integrator as the underlying numerical scheme, we ensure that the resulting uncertainty propagation method inherits the geometric properties of the time evolution of probability densities in Hamiltonian systems, that arise from the symplectic geometry of the phase space.

A natural application of the proposed scheme is to the problem of global attitude estimation, particularly when the dynamics of the rigid body are extremely nonlinear, and the attitude measurements are relatively infrequent. Estimation problems with such characteristics are problematic for traditional techniques, such as the Kalman filter, which require frequent measurements, and relatively benign dynamics, in order to justify the localization and linearization assumptions built into the method.
\bibliography{est}
\bibliographystyle{IEEEtran}

\end{document}